\documentclass[letter, 10pt, conference]{ieeeconf}  
\IEEEoverridecommandlockouts
\usepackage{makecell}

\usepackage{amsmath,amssymb,amsfonts}
\usepackage{algorithmic}
\usepackage{graphicx}
\usepackage{textcomp}
\usepackage{subcaption}
\usepackage{xcolor}
\usepackage{algorithm} 
\usepackage{amsmath}
\usepackage{cancel}
\usepackage{mathrsfs}
\usepackage{cite}



\usepackage{amsthm}

\def\idty{{\operatorname{Id}}}

\def\cF{{\mathcal F}}
\def\cZ{{\mathcal Z}}

\def\R{{\mathbb R}}
\def\S{{\mathbb S}}

  \def\cG{{\mathcal G}}  \def\cS{{\mathcal S}}    \def\cN{{\mathcal N}}                 \def\cF{{\mathcal F}}    \def\cX{{\mathcal X}} \def\cY{{\mathcal Y}}  \def\cZ{{\mathcal Z}}

\theoremstyle{plain}
\newtheorem{theorem}{Theorem}[section] 

\theoremstyle{definition}
\newtheorem{definition}{Definition}
\newtheorem{notations}[definition]{Notation}
\newtheorem{remark}[theorem]{Remark}
\newtheorem{assumption}[theorem]{Assumptions}

\title{\LARGE \bf
Scalable iterative Gramian synthesis for control-affine systems
}

\author{Zongxi Yu, Cyprien Tamekue, Ruiqi Chen and ShiNung Ching
\thanks{This work has been submitted to the IEEE for possible publication. Copyright may be transferred without notice, after which this version may no longer be accessible.}
\thanks{This work is partially supported by grant R21MH132240 from the US National Institutes of Health to SC.}
\thanks{Zongxi Yu, Cyprien Tamekue, and ShiNung Ching are with the Department of Electrical and Systems Engineering, Washington University in St. Louis, St. Louis, 63130, MO, USA; Ruiqi Chen is with the Neurosciences Program in the Division of Biology and Biomedical Sciences, Washington University in St. Louis, St. Louis, 63130, MO, USA;
        {\tt\small e-mail: \{y.zongxi,cyprien,chen.ruiqi,shinung\}@wustl.edu.}}%
}

\begin{document}

\maketitle

\begin{abstract}
This article presents a scalable implementation of nonlinear Gramian-based control synthesis for control-affine systems, including a minimum energy control construction. These synthesis advances are achieved by addressing key computational bottlenecks inherent to iterative synthesis map formulations, yielding 
a computational scheme that exhibits rapid convergence and high-precision.
The efficacy of this synthesis framework is demonstrated across five canonical nonlinear control systems and 100-dimensional recurrent neural network models, including underactuated systems. Empirical scaling results further indicate that convergence is primarily governed by intrinsic system properties, such as nonlinearity and controllability, rather than by state-space dimensionality. This work provides a practical, scalable computational pathway for translating rigorous nonlinear synthesis theory into high-dimensional control applications.
\end{abstract}

\section{INTRODUCTION}\label{s:introduction}


We consider
control-affine systems 
\begin{equation}\label{eq:state-dependent-input}
\dot{x}(t)=N_t(x(t))+B_t(x(t))\,u(t),\qquad t\in[t_0,T],
\end{equation}
where $0\le t_0< T$, $x(t)\in\mathbb R^d$ is the state, $u(t)\in\mathbb R^k$ the control input with $k\le d$, $N_t$ is a possibly nonautonomous drift and $B_t$ is a possibly time- and state-dependent input matrix. Systems of the form \eqref{eq:state-dependent-input} arise throughout engineering and applied science~\cite{cheng2010analysis,coron2007control}.
They also encapsulate a wide range of emerging applications in neural networks and neural engineering, in which the state vector corresponds, e.g., to a set of neurons or brain regions (see, e.g., \cite{chenDynamicalModelsReveal2025a,singh2020estimation}). In such applicative contexts, the ability to generate controls must scale to dimensionality in the hundreds or greater.




Control analysis and design of nonlinear control-affine systems has been recognized as theoretically nontrivial. In the linear time-invariant (LTI) and time-variant (LTV) settings, controllability admits algebraic characterizations and constructive synthesis via the controllability Gramian~\cite{kalman1963mathematical}.
By contrast, for nonlinear control-affine systems, control analysis and synthesis depend on the state, the input, and the underlying trajectory, so constructive synthesis is typically more delicate and often relies on geometric arguments, continuation methods, or problem-specific optimal-control routines~\cite{agrachev2013geometric,chitour2006continuation,hermann2003nonlinear,sun2007necessary,sussmann1992new,ji2023global,jean2014nonholonomic,trelat2012optimal,zuyev2017motion}.

Our recent work~\cite{tamekue2025controlanalysis,tamekue2026minimum} has begun to recover a nonlinear analogue of Gramian calculus for~\eqref{eq:state-dependent-input} by introducing trajectory-dependent Gramian maps and casting steering as a fixed-point problem on feasible coercivity classes---subsets on which the Gramian is uniformly coercive/invertible. In that framework, one obtains a constructive Gramian-like controller with an intrinsic energy certificate and a Picard convergence mechanism.
A complementary development~\cite{tamekue2026minimum} shows that minimum-energy controls can also be characterized as the fixed-point of a Lagrange-multiplier map involving a generally non-symmetric trajectory-dependent Gramian, and that this minimum-energy controller coincides with the symmetric Gramian controller~\cite{tamekue2025controlanalysis} under an orthogonality condition.
In this work, we refer to the first framework~\cite{tamekue2025controlanalysis} as General Synthesis, and the latest development~\cite{tamekue2026minimum}  as Minimum-energy Synthesis.

However, translating this Gramian synthesis framework into a scalable numerical method---one applicable to neural models, for instance---is numerically nontrivial. First, the synthesized control is defined pointwise via the Jacobian flow evaluated along the steered trajectory and lives in an infinite-dimensional function space, making its numerical representation and evaluation computationally costly. Second, the nonlinear controllability Gramian integrates the Jacobian flow along the steered trajectory and admits no closed-form ODE; it must be approximated numerically via quadrature. Thus, a precise approximation is computationally costly in high dimensions, especially over a long time horizon. 


In this work, we present a scalable implementation of the general Gramian synthesis framework \cite{tamekue2025controlanalysis} and provide the first numerical implementation of the minimum-energy Gramian synthesis theory \cite{tamekue2026minimum}. We validate both syntheses across canonical nonlinear control examples, and demonstrate scalability via controlling a 100-dimensional {data-driven, and possibly underactuated} brain network model \cite{chenDynamicalModelsReveal2025a}. 
This work demonstrates the potential of modernizing traditional numerical simulation workflows with {sophisticated algorithmic design, as well as adopting} state-of-the-art computational {infrastructure}, enabling {more rigorous control theories like the Gramian Synthesis framework to reach broader applied and interdisciplinary audiences.




\section{Preliminaries on nonlinear Gramian synthesis }
\label{s:preliminary}
This section briefly recalls the two nonlinear Gramian-based fixed-point controllers that motivate this work. Our purpose here is not to restate the full existence theory~\cite{tamekue2025controlanalysis,tamekue2026minimum}, but to isolate the common computational structure behind the two synthesis maps. Throughout, we consider the control-affine system~\eqref{eq:state-dependent-input} with boundary values
\begin{equation}\label{eq:endpoints}
    x(t_0)=x^0\in\R^d,\qquad x(T) = x^1\in\R^d.
\end{equation}

\begin{assumption}\label{ass:on N_t}
 The nonautonomous drift $N_t:[t_0, T]\times\R^d\to\R^d$ satisfies the following regularities assumptions
 \begin{enumerate}
        \item The map $t\mapsto N_t(x)$ is $L^\infty$, for every fixed $x\in\R^d$,
        \item The map $x\mapsto N_t(x)$ is $C^2$, for every fixed $t\in[t_0, T]$.
    \end{enumerate}
    Additionally, for some constants $\Lambda_1\ge0$ and $\Lambda_2\ge0$,
    \begin{equation}\label{eq:estimates on N_t}
        \|D_xN_t(x)\|\le\Lambda_1,\qquad\|D_x^2N_t(x)\|\le\Lambda_2,
    \end{equation}
for all $(t,x)\in[t_0, T]\times\R^d$.
    Here, $D_xN_t(x)$ is the differential of $N_t$, and $D_x^2N_t(x)$ its second-derivative at any $x\in\R^d$.
\end{assumption}

\begin{assumption}\label{ass:on B_t}
     The input matrix $B_t\in L^\infty((t_0, T)\times\R^d; \R^{d\times k})$. Furthermore, for a fixed $t\in[t_0, T]$, $x\mapsto B_t(x)$ is $C^2$, and for some constants $L_B\ge0$ and $L_{DB}\ge0$, it holds for all $(t,x)\in[t_0, T]\times\R^d$,
\begin{equation}\label{eq:estimates on B_t}
    \|D_xB_t(x)\|\le L_B,\qquad \|D_x^2B_t(x)\|\le L_{DB}.
\end{equation}
\end{assumption}

\begin{remark}
    The uniform bounds~\eqref{eq:estimates on N_t} and~\eqref{eq:estimates on B_t} are mainly a technical convenience, and all derivative bounds are to be understood either globally on $\R^d$ or locally on a compact set containing the trajectories of interest over $[t_0, T]$. 
\end{remark}


\begin{notations}
    Throughout the following, we fix $i\in\{1,2\}$, $\tau_1=t_0$ and $\tau_2=T$. We also fix 
    $\cX:=L^2((t_0, T); \R^k)$, $\cY:=L^\infty((t_0, T); \R^k)$. 
    Finally, we use $|x|$ for the Euclidean norm of $x\in\R^d$, and  $\mathscr{L}(X, Y)$ as the space of linear bounded operators from two normed vector spaces $X$ and $Y$. 
\end{notations}

\subsection{Flow-conjugate solution representation}

Let $x\mapsto\Phi_{s,t}(x)$ denote the flow of the drift equation,
\begin{equation}\label{eq:nonautonomous flow}
    \partial_t\Phi_{s,t}(x^0) = N_t(\Phi_{s,t}(x)),\quad \Phi_{s,s}(x) = x\in\R^d
\end{equation}
where $(s,t)\in[t_0,T]^2$. A key starting point in the nonlinear Gramian framework~\cite{tamekue2025controlanalysis,tamekue2026minimum} is the flow-conjugate representation of trajectories of~\eqref{eq:state-dependent-input}. For $x^0\in\R^d$ and control $u\in\cY$,
the corresponding trajectory $x_u(\cdot)\in C^0([t_0,T]; \R^d)$ to~\eqref{eq:state-dependent-input} can be written in the form (see~\cite[Theorem 3.2]{tamekue2025controlanalysis})
 \begin{equation}\label{eq:flow_conjugate}
    \begin{gathered}
         x_u(t)=\Phi_{\tau_i,t}\left(\Phi_{t_0,\tau_i}(x^0)+I_u(t)\right),\quad\forall t\in[t_0, T],\\
        I_u(t):=\int_{t_0}^{t}\!\!\! D\Phi_{s,\tau_i}\big(x_u(s)\big)B_s(x_u(s))u(s)\,ds
    \end{gathered}
    \end{equation}
    where $D\Phi_{s,\tau_i}\big(x_u(s)\big)$ is the differential of $\Phi_{t,\tau_i}$ at $x_u(s)$.

For ease in notation, we define the residual vectors
\begin{equation}\label{eq:yi_def}
y_i:=\Phi_{T,\tau_i}(x^1)-\Phi_{t_0,\tau_i}(x^0),\qquad i\in\{1,2\}.
\end{equation}
Then, steering~\eqref{eq:state-dependent-input} from $x^0$ to $x^1$ over $[t_0, T]$ is equivalent to solving an implicit integral equation of the form
\begin{equation}\label{eq:Lu_tau}
L_{u,\tau_i}u=y_i,
\end{equation}
where $L_{u,\tau_i}\in\mathscr{L}(\cX,\R^d)$ is defined by
\begin{equation}\label{eq:map L_u tau_i}
    \hspace{-0.1cm}L_{u,\tau_i}v = \int_{t_0}^T\!\!\!D\Phi_{t,\tau_i}(x_{u}(t))\,B_t(x_u(t))\,v(t)\,dt,\quad v\in\cX.
\end{equation}

\subsection{Trajectory-dependent Gramian maps}

For a given $x^0\in\R^d$ and control $u\in\cY$, let $x_u(\cdot)\in C^0([t_0,T]; \R^d)$ denote the corresponding trajectory of~\eqref{eq:state-dependent-input} given by~\eqref{eq:flow_conjugate}. 
The first, symmetric, Gramian-like matrix is
\begin{equation}\label{eq:almost-optimal Gramian}
   \cN_{\tau_i}(u) = L_{u,\tau_i}L_{u,\tau_i}^\ast
\end{equation}
where $L_{u,\tau_i}^\ast\in\mathscr{L}(\R^d,\cX)$ denote the adjoint of $L_{u,\tau_i}$, namely
\begin{equation}
    [L_{u,\tau_i}^\ast z](t) =  B_t(x_u(t))^\top D\Phi_{t,\tau_i}(x_{u}(t))^\top z,\;\, z\in\R^d.
\end{equation}
The matrix $\cN_{\tau_i}(u)\in\R^{d\times d}$ is symmetric positive semidefinite and is the nonlinear analogue of the controllability Gramian~\cite{kalman1963mathematical} used in the constructive steering map of~\cite{tamekue2025controlanalysis}.

The minimum-energy synthesis~\cite{tamekue2026minimum} introduces, in addition, a mixed, generally non-symmetric Gramian
\begin{equation}\label{eq:optimal Gramian}
     \cG_{\tau_i}(u) = L_{u,\tau_i}DF_{\tau_i}(u)^\ast
\end{equation}
where $F_{\tau_i}:u\in\cY\subset\cX\mapsto\R^d$ is the map defined by
\begin{equation}\label{eq:feasible map}
    F_{\tau_i}(u) = L_{u,\tau_i}u-y_i
\end{equation}
and $\,DF_{\tau_i}(u)^\ast\in\mathscr{L}(\R^d,\cX)$ denote the adjoint of the differential $DF_{\tau_i}(u)\in\mathscr{L}(\cX,\R^d)$ of $F_{\tau_i}$ at $u\in\cY$, namely
\begin{equation}
\label{eq:DF}
    \begin{gathered}
        [DF_{\tau_i}(u)^\ast z](t) =  B_t(x_u(t))^\top Q_{u,\tau_i}(T,t)^\top z\\
        Q_{u,\tau_i}(T,t):=D\Phi_{T,\tau_i}(x_u(T))R_u(T,t)
    \end{gathered}
\end{equation}
with $R_u(t,s)\in C^0([t_0, T]^2;\R^{d\times d})$ is the state-transition matrix (STM) of the linearized equation
 \begin{equation*}\label{eq:linearized equation B=0}
       \dot{y}(t) = D_x\big[N_t(x_u(t))+B_t(x_u(t))u(t)\big] y(t),\;t\in[t_0, T]
    \end{equation*}
    satisfying $R_u(s,s)=\idty$.


\subsection{The General Synthesis}\label{ss:symmetric Gramian-like}

The first synthesis map, introduced in~\cite{tamekue2025controlanalysis}, is 
\begin{equation}\label{eq:Si_map}
\cS_i(u):=L_{u,\tau_i}^\ast\,\cN_{\tau_i}(u)^{-1}y_i,\qquad \forall u\in\cF_i\subset\cY
\end{equation}
where $\cF_i\subset\cY$ is a subset on which $\cN_i$ is uniformly coercive. The fixed point $u_i=\cS_i(u_i)$ yields a steering control for \eqref{eq:state-dependent-input}, together with an intrinsic energy certificate~\cite[Theorem 3.15]{tamekue2025controlanalysis} whenever $\cF_i\neq\emptyset$ and $\cS_i$ is a self-map on $\cF_i$. Moreover, $u_i$ is obtained as the limit of the Picard iteration $u^{(n+1)}=\cS_i(u^{(n)})$ for any $u^{(0)}\in\cF_i$.

For computation, $u^{(n+1)}=\cS_i(u^{(n)})$ should be read as follows: given a current control iterate $u:=u^{(n)}$,

\noindent1. solve the state equation~\eqref{eq:state-dependent-input} to obtain $x_u(\cdot)$;

\noindent 2. propagate the associated variational equation
\begin{equation}\label{eq:variational propagation}
    \begin{cases}
        \dot y(s) = N_s(y(s)),\quad y(t) = x_u(t)\\
        \dot Y(s) = D_xN_s(y(s))Y(s),\quad Y(t) = B_t(x_u(t))
    \end{cases}
\end{equation}
forward on $[t, T]$ or backward on $[t_0, t]$ to compute for a fixed $t\in [t_0, T]$, the product $D\Phi_{t,T}(x_{u}(t))\,B_t(x_u(t))$ or $D\Phi_{t,t_0}(x_{u}(t))\,B_t(x_u(t))$;

\noindent3. assemble the symmetric Gramian $\cN_{\tau_i}(u)$;

\noindent4. solve the linear system $\cN_{\tau_i}(u)\,\lambda = y_i$;

\noindent5. form the updated control $u^{(n+1)}=[L_{u,\tau_i}^\ast\lambda](t)$.

This five-step procedure is summarized in Algorithm~\ref{alg:picard}.

\subsection{The Minimum-energy Synthesis}\label{ss:non-symmetric Gramian-like}

The second synthesis, derived in~\cite{tamekue2026minimum}, 
is built on the Lagrange multiplier map
\begin{equation}\label{eq:Zi_map}
\cZ_i(u):=DF_{\tau_i}(u)^\ast\,\cG_{\tau_i}(u)^{-1}y_i,\qquad\forall u\in\cF_i.
\end{equation}
Here, $\cF_i\subset\cY$ is identical to~\eqref{eq:Si_map} and corresponds to a unified subset on which $\cN_{\tau_i}$ (resp. $\cG_{\tau_i}$)  is uniformly coercive (resp. invertible). When $\cF_i\neq\emptyset$ and $\cZ_i$ is a self-map on $\cF_i$, $\cZ_i$ admits a unique fixed point $\bar u_i$, and the corresponding Picard iteration $u^{(n+1)}=\cZ_i(u^{(n)})$, $u^{(0)}\in\cF_i$ converges to $\bar u_i$ with factorial decay. 
Algorithmically, $u^{(n+1)}=\cZ_i(u^{(n)})$ follows the same computational pattern as previously: given a current control iterate $u:=u^{(n)}$,

\noindent1. solve the state equation~\eqref{eq:state-dependent-input} to obtain $x_u(\cdot)$;

\noindent 2. propagate the variational equation~\eqref{eq:variational propagation} to compute $D\Phi_{t,\tau_i}(x_{u}(t))\,B_t(x_u(t))$, the STM-related equation
\begin{equation}\label{eq:variational propagation STM}
    \begin{cases}
        \dot Y(s) = D_x\big[N_s(x_u(s))+B_s(x_u(s))u(s)\big]Y(s),\\
         Y(t) = B_t(x_u(t))
    \end{cases}
\end{equation}
forward on $[t, T]$ to compute for a fixed $t\in [t_0, T]$, the product $\Lambda_u(t,T):=R_u(T,t)\,B_t(x_u(t))$, or
\begin{equation}\label{eq:variational propagation backward}
    \begin{cases}
        \dot y(s) = N_s(y(s)),\quad y(T) = x_u(T)\\
        \dot Y_t(s) = D_xN_s(y(s))Y_t(s),\quad Y_t(T) = \Lambda_u(t,T)
    \end{cases}
\end{equation}
backward on $[t_0,T]$ to compute for a fixed $t\in [t_0, T]$,  the product $D\Phi_{T,t_0}(x_u(T))R_u(T,t)\,B_t(x_u(t))$;

\noindent3. assemble the non-symmetric Gramian $\cG_{\tau_i}(u)$;

\noindent4. solve the linear system $\cG_{\tau_i}(u)\,\lambda = y_i$;

\noindent5. form the updated control $u^{(n+1)}(t)=[DF_{\tau_i}(u)^\ast\lambda](t)$.

This five-step procedure is also summarized in Algorithm~\ref{alg:picard}.

\section{Implementation of Nonlinear Gramian Synthesis}

The Picard algorithms and implementations presented in this section assume that for the given boundary condition~\eqref{eq:endpoints}, there exists a unified nonempty subset $\cF_i\subset\cY$ over which the synthesis map $\cS_i$ and $\cZ_i$ are self-maps so that 
from any initial reference control $u_i^{(0)}\in\cF_i$, the convergence of the Picard iterate to the fixed point is guaranteed. 

\subsection{Picard iteration procedure}
Algorithm \ref{alg:picard} summarizes the Picard Gramian Synthesis procedure. During the iterations, we record the terminal-state (end-point) error.
\[
\mathrm{err}_{\rm end}(u^{(n)}):=|x_{u^{(n)}}(T)-x^1|,
\]
the control update
\[
\mathrm{err}_{\rm fp}(u^{(n)}):=\operatorname*{ess\,sup}_{t\in[t_0,T]}
|u^{(n+1)}(t)-u^{(n)}(t)|,
\]
and the control energy
\[
E(u^{(n)})=\frac{1}{2}\int_{t_0}^T |u^{(n)}(t)|^2\,dt.
\]

\begin{algorithm}
\caption{Picard iteration for $\mathcal S_i$ (General Synthesis) and $\mathcal Z_i$ (Minimum-energy Synthesis)}
\label{alg:picard}
\begin{algorithmic}[1]

\REQUIRE Boundary condition $x^0$, $x^1$, $t_0$, $T$

\REQUIRE$u^{(0)}(\cdot)$, $N_{\max}$, $K$, $\epsilon$, $i\in \{1,2\}$ , $\mathrm{Map}\in\{\cS_i,\cZ_i\}$
\ENSURE $\cN_{\tau_i}(u^{(0)})$ or $\cG_{\tau_i}(u^{(0)})$ is invertible


\IF{$i=1$}
    \STATE $\tau_i \gets t_0$
\ELSIF{$i=2$ }
 \STATE $\tau_i \gets T$
 \ENDIF
\STATE $y_i \leftarrow \Phi_{T,\tau_i}(x^1) - \Phi_{t_0,\tau_i}(x^0)$ as \eqref{eq:yi_def} 
 \STATE Let $\{t_k\}_{k=1}^{K}$ be quadrature points over $[t_0, T]$

\FOR{$n=0,\dots,N_{\max}-1$}
    \STATE $u(\cdot) \gets u^{(n)}(\cdot)$    
   
    \FOR{$k = 1, \dots, K$} 

    \STATE Solve $x_u(t)$ from $t_0$ to $t_k$ as ~\eqref{eq:state-dependent-input}
    \STATE Solve  $D\Phi_{t_k,\tau_i}(x_{u}(t_k))\,B_{t_k} (x_u(t_k))$ from $t_k$ to $\tau_i$ as ~\eqref{eq:variational propagation}     \textit{// shared by both syntheses}
    \ENDFOR
    \IF{$\mathrm{Map}=\mathcal S_i$}

        \STATE Assemble $\mathcal N_{\tau_i}(u)$ as~\eqref{eq:almost-optimal Gramian} over the $K$ samples

        \STATE Solve $\cN_{\tau_i}(u)\,\lambda = y_i$ 
        \STATE $u^{(n+1)}(\cdot)\gets B_t(x_u(\cdot))^\top\,D\Phi_{t,\tau_i}(x_{u}(\cdot))^\top\,\lambda$ 
    \ELSIF{$\mathrm{Map}=\mathcal Z_i$}
     \FOR{$k = 1, \dots, K$}
        \STATE Solve $D\Phi_{T,\tau_i}(x_u(T))R_u(T,t_k)B_t(x_u(t_k))$ from $t_k$ to $\tau_i$ as~\eqref{eq:variational propagation STM} or~\eqref{eq:variational propagation backward}
    \ENDFOR
        \STATE Assemble $\cG_{\tau_i}(u)$ as~\eqref{eq:optimal Gramian} over the $K$ samples
        \STATE Solve $\cG_{\tau_i}(u)\,\lambda = y_i$ 
        \STATE $u^{(n+1)}(\cdot)\gets B_t(x_u(\cdot))^\top\,Q_{u,\tau_i}(T,\cdot)^\top\,\lambda$ 
    \ENDIF

    \IF{termination criterion holds}
        \STATE \textbf{break}
    \ENDIF
\ENDFOR
\end{algorithmic}
\end{algorithm}

\textbf{Termination criteria:}
The iteration terminates when any of the following conditions is met:
\begin{enumerate}
    \item Maximum iteration number is reached: $n \ge N_{\max}$;
    \item End-point error converges: $\mathrm{err}_{\rm end}(u^{(n)}) \le \epsilon_x$; 
    \item Control function converges: $\mathrm{err}_{\rm fp}(u^{(n)}) \le \epsilon_u$.
\end{enumerate}
Here, $N_{\max}$ is the maximum number of iterations, $\epsilon_x>0$ is the terminal-state tolerance, and $\epsilon_u>0$ is the control-update tolerance. 

\subsection{Numerical implementation}
We implement the Gramian synthesis using JAX~\cite{jax2018github} and Diffrax~\cite{kidger2021on}. We use the 8th-order Dormand-Prince Runge-Kutta method from Diffrax as the ODE solver, with step size controlled by a PID controller with default parameters $(p,i,d)=(0,1,0)$. The numerical quadrature uses a Simpson's rule, as implemented in the JAX-native Quadax library~\cite{conlin_2026_19059171}. Matrix inversion uses Cholesky decomposition for well-conditioned, positive-semidefinite Gramian matrices and a least-squares solver otherwise. 

In the rest of this section, we present several key technical design choices in our implementation, which are generalizable to other nonlinear control simulations. In particular, our techniques leverage modern computational tools, including automatic differentiation, vectorization, GPU acceleration, and kernel fusion, to address numerical constraints arising from the mathematical structure of Gramian synthesis. 

\paragraph{Auto-differentiation instead of analytical gradient} In our implementation, Jacobians of the drift $D_xN_s(x)$ and controlled trajectory $D_x\big[N_t(x)+B_t(x)u(s)\big]$ are computed via forward-mode automatic differentiation.

\paragraph{Augmented variational equation over reverse-mode automatic differentiation
}
Although the flow Jacobian $D\Phi_{t,\tau_i}(x_{u}(t))$ can be obtained by automatic differentiation through the ODE solver~\cite{kidger2021on}, these methods either store intermediate states at a memory cost of $\mathcal{O}(dR+d^2)$, where $R$ is the total number of solver steps,  or introduce additional numerical error via backward integration. Instead, we compute the flow Jacobian analytically by directly solving variational equations, thereby maintaining an overall memory cost of $\mathcal{O}(d^2)$ and significantly reducing peak memory usage for high-dimensional or highly nonlinear systems.

\paragraph{Analytical fusion of mathematical operations}
Since the flow Jacobian $D\Phi_{t,\tau_i}(x_{u}(t))$ always appears in the product form $D\Phi_{t,\tau_i}(x_u(t))B(x_u(t))$, we solve directly for this product via an augmented ODE \eqref{eq:variational propagation}, reducing memory complexity from $\mathcal{O}(d^2)$ to $\mathcal{O}(dk)$, which is significant for underactuated systems where $k\ll d$. The flow Jacobian computation shares an identical computational graph, differing only in the initial state and time horizon, and is invariant across Picard iterates with different control $u(\cdot)$. Thus, we use just-in-time (JIT)~\cite{jax2018github} compilation to speed up repeated calls, which further benefits from the reduced dimensionality of the product form through more aggressive kernel fusion.

\paragraph{Parallelization of quadrature integration via vectorization}
In each Picard iteration, the controllability Gramian is approximated via numerical quadrature over $K$ sample points, where each sample requires evaluating the Jacobian product $D\Phi_{t,\tau_i}(x_u(t))B(x_u(t))$ via an ODE solve. A large $K$ is unavoidable for nonlinear systems, as the integrand varies significantly across the time horizon, with approximation error propagating quickly with each Gramian inversion. Rather than evaluating samples sequentially, we vectorize the entire quadrature computation via \texttt{jax.vmap}, parallelizing all $K$ ODE solves simultaneously. Such vectorization decouples quadrature accuracy from wall-clock time, allowing arbitrarily large $K$ without proportional increase 
in runtime within hardware constraints, and potentially enabling future horizontal scaling across 
multiple accelerators via JAX's built-in parallel/distributed computing features.

\paragraph{Efficient and precise infinite-dimensional function representation}

The total number of solver calls of a precise yet naive implementation of Picard iteration grows exponentially with the iteration number, since $u^{(n)}(t)$ depends on $x_u^{(n-1)}(t)$, which itself depends on $u^{(n-2)}(t)$ recursively; each evaluation triggers a chain of nested ODE solves. We break the recursive dependence by representing the controlled trajectory $x_u^{(n)}(t)$ via the dense output interpolant produced by the adaptive solver, decoupling each Picard iteration from all prior ODE solves while retaining the target precision. This change allows us to maintain a similar algorithmic cost per iteration under a stable solver step size.

Traditionally, a continuous-time control function is approximated by a finite basis or local interpolation over a predefined grid sequentially, introducing a tradeoff between approximation error and computational cost that cannot be resolved a priori. We instead offer two evaluation strategies, selected based on the system's nonlinearity, dimensionality, and target precision. 

\textbf{On-demand evaluation} The control $u^{(n)}(t)$ is lazily evaluated at each solver-selected time point, delegating temporal resolution entirely to the adaptive step-size controller. This aligns the control's approximation accuracy with the solver's error tolerance. However, such evaluation would occur sequentially, which is not suitable for situations using large solver steps.

\textbf{Dense interpolation} For systems with strong nonlinearity where on-demand evaluation becomes intractable, especially since each evaluation of $u(t)$ requires solving $D\Phi_{t,\tau_i}(x_u(t))B(x_u(t))$ or $R_u(T,t)$ (\eqref{eq:DF}) with different initial conditions. We precompute $u(t)$ on a grid and store a dense interpolant, vectorized via \texttt{jax.vmap}, conceptually similar to batching or caching. This technique trades some approximation accuracy for an $\mathcal{O}(\log{M})$ evaluation cost and decouples $K$, the interpolation sample size, from the actual runtime.

\section{Simulations and results}\label{s:simulations and results}
\begin{table}[t]
\centering
\renewcommand{\arraystretch}{1.3}
\resizebox{\columnwidth}{!}{%
\begin{tabular}{lccccc}
\hline
System & $d$ & $k$ & $x^0$ & $x^1$ & $[t_0, T]$ \\
\hline
Pendulum 
    & 2 & 1 
    & $[0,\, 0]^\top$ 
    & $[\pi,\, 0]^\top$ 
    & $[0.5,\, 1.5]$ \\
Unicycle 
    & 3 & 2 
    & $[0.5,\, 0.25,\, \frac{\pi}{12}]^\top$ 
    & $[1,\, 0.75,\, \frac{4\pi}{3}]^\top$ 
    & $[0,\, 2]$ \\
Spacecraft
    & 6 & 3 
    & $[0.3,\,0.2,\,0.1,\,0,\,0,\,0]^\top$ 
    & $[0,\,0,\,0,\,0,\,0,\,0]^\top$ 
    & $[0,\, 5]$ \\
SIR
    & 3 & 1 
    & $[1,\,0.2,\,0.1]^\top$ 
    & $[0.5,\,0.25,\,0.2]^\top$ 
    & $[0,\, 0.5]$ \\

2D Hopfield
& 2 & 1, 2 
& $[1.0,1.0]^\top$ 
& $[-1.0,\, -1.0]^\top$ 
& $[0,\, 1.5]$  \\
MINDy 
    & 100 & 100 
    & $\mathbf{1}_{100}$ 
    & $-0.5\cdot\mathbf{1}_{100}$
    & $[0,\,3]$ \\
MINDy
    & 100 & 50
    & $\sim \mathcal{N}(\mathbf{0}_{100}, \mathbf{I}_{100})$
    & $x_{u_{\text{sim}}}(T)$ (see Section~\ref{sec:underactuated})
    & $[0,\,4]$ \\
\hline
\end{tabular}}
\caption{Simulation setup for benchmark systems.}
\label{tab:setup}
\end{table}
\subsection{System setup}

We evaluate the nonlinear Gramian synthesis across six example systems spanning a range of structural properties: 
time-varying dynamics, state-dependent input matrices, driftless 
dynamics, oscillatory behavior, and high dimensionality. 
Together, the systems provide coverage of non-autonomous and 
state-dependent control-affine structures, with applications in 
motion planning, epidemiology, and neuroscience. 
\paragraph{Unicycle}
The kinematic unicycle model is a driftless control-affine system
with state \(x:=(p,\theta)\in\R^2\times\S^1\) and controlled via $u:=(v,\omega)^\top\in\R^2$, namely
\begin{equation}
N_t(x)=0_{\R^3} 
\qquad
B_t(x)^\top=\begin{bmatrix}\cos(\theta)& \sin(\theta)&0\\0 & 0 & 1\end{bmatrix}.
\end{equation}
\paragraph{Torque-controlled pendulum with varying length}
We consider the nondimensional equation of a torque-controlled pendulum with varying length~\cite[Section~2]{belyakov2009dynamics},
\begin{equation}
N_t(x)=\begin{bmatrix}x_2\\-a(t)\sin x_1-\gamma(t)x_2\end{bmatrix},
\;\,
B_t(x)=\begin{bmatrix}0\\b(t)\end{bmatrix},
\end{equation}
where
$
b(t)=\bigl(1+0.5\cos(t)\bigr)^{-2}$, $a(t)=\lambda^2\sqrt{b(t)}$,
$\gamma(t)=- b(t)\sin(t)+\beta\lambda$ with $\lambda=0.78$ and $\beta=0.13$.



\paragraph{SIR Epidemic Model}
We consider a SIR epidemic model in $3$D with state $x:=(S,I,R)\in\R^3$, and a scalar control $u\in\R$ that directly acts on the $S$ compartment, 
\begin{equation}
N_t(x)=\begin{bmatrix}
    \lambda - \beta S I - \mu S\\
    \beta S I - (\mu+\gamma)I\\
    \gamma I - \mu R
\end{bmatrix},
\qquad
B_t(x)=\begin{bmatrix}-S\\0\\0\end{bmatrix},
\end{equation}
where $\lambda=1$, $\beta=2$, $\mu=0.2$ and $\gamma=1$.
\paragraph{Spacecraft Attitude Control}
We consider the control of the attitude of a rigid
spacecraft with control torques provided by thruster jets, as described in~\cite[Example~3.9]{coron2007control}. The state $x:=(\eta,\omega)\in\R^3\times\R^3$,
where $\eta=(\phi,\theta,\psi)^\top$ and
$\omega=(w_1,w_2,w_3)^\top$, and the control $u\in\R^3$ with $B_t(x)^\top = \begin{bmatrix}0_{3\times 3} & J^{-1}\end{bmatrix}$ and $J = \operatorname{diag}(10,20,15)$.

\paragraph{2D Hopfield system} We consider the following fully- or under-actuacted 2D Hopfield neural network~\cite{hopfield1982neural}, 
\begin{equation}
\begin{aligned}
N_t(x)&=-\begin{bmatrix}
    0.5x_1\\0.3x_2
\end{bmatrix}+\begin{bmatrix}
    0.5&-1.5\\
    1.5&-0.5
\end{bmatrix}\begin{bmatrix}
    \tanh(x_1)\\
    \tanh(x_2)
\end{bmatrix},\\
B_t(x)&=\begin{bmatrix}
    1 & 0.5
\end{bmatrix}^\top\qquad\text{or}\qquad B_t(x)= \idty.
\end{aligned}
\end{equation}

\paragraph{High-dimensional MINDy fMRI model}
The MINDy model is a Hopfield-type neural network~\cite{hopfield1982neural} with parameters estimated from resting-state fMRI signals~\cite{chenDynamicalModelsReveal2025a}, namely,
\begin{equation}
N_t(x)=-D x+\mathbf{W}\psi_{\alpha,\beta}(x),
\qquad
B_t(x)=\idty
\end{equation}
where $D\in\mathbb{R}^{d\times d}$ is the diagonal decay matrix, $\mathbf{W}\in\mathbb{R}^{d\times d}$ is a synaptic weight matrix, and $\psi_{\alpha,\beta}(x)=\sqrt{\alpha^2+(\beta x+0.5)^2}-\sqrt{\alpha^2+(\beta x-0.5)^2}$ is an elementwise activation parameterized with $\alpha\in\R^d$ and $\beta=20/3$.
\subsection{Picard convergence across benchmark systems}

\begin{figure*}
    \centering
    \includegraphics[width=1.0\linewidth]{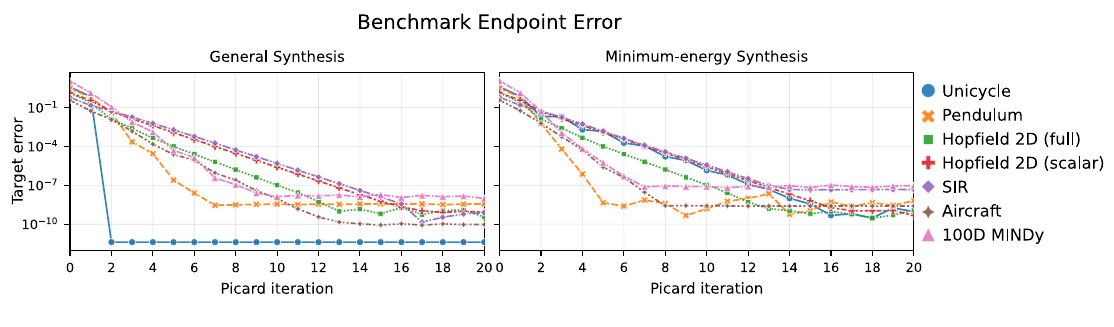}
    \caption{End-point error convergence for benchmark systems}
    \label{fig:convergence}
\end{figure*}

We simulated both general and minimum-energy Gramian synthesis with the boundary conditions summarized in Table~\ref{tab:setup}, using 20 fixed Picard iterations per system. Fig.~\ref{fig:convergence} reports the endpoint error $|x_u^{(n)}(t_1) - x^1|$ across iterations. To ensure computational tractability, we use double precision with $\text{rtol}=10^{-7}$, $\text{atol}=10^{-7}$ for the minimum-energy synthesis on the 100D MINDy model, and double precision with $\text{rtol}=10^{-8}$, $\text{atol}=10^{-10}$ for all other simulations. We observed convergence consistently across all cases, with the general Gramian synthesis reaching endpoint errors of $\sim 10^{-9}$--$10^{-11}$, approaching the expected numerical precision floor. The minimum-energy Gramian converges to a higher yet still acceptable endpoint error. In both cases, the convergence curves are consistent with the factorial decay guaranteed by the theory~\cite{tamekue2025controlanalysis, tamekue2026minimum}. {To
isolate the effect of actuation, we evaluated two configurations
of the 2-dimensional Hopfield network with identical dynamics: an underactuated configuration with scalar input  and a fully actuated configuration with two inputs. The fully actuated configuration converges faster than the under-actuated counterpart in both synthesis, substantiating the controllability effect on convergence. }

A key empirical finding is that, for the syntheses, convergence is primarily governed by the nonlinearity of the vector field and the system's intrinsic controllability, rather than by state-space dimension. The unicycle converges within 2 iterations for general Gramian synthesis, reflecting its driftless structure and full controllability. The 100D MINDy model converges comparably to the torque-controlled pendulum, both reaching their precision floors within iteration 10. The SIR and spacecraft models require the most iterations ($\sim$15), suggesting that the numerical conditions of the intrinsic controllability govern convergence more than dimensionality.


Taken together, our empirical results confirm that both syntheses are computable across all benchmark systems and consistent with their respective theoretical guarantees~\cite{tamekue2025controlanalysis, tamekue2026minimum}, representing complementary computational tools for the nonlinear control community.
\subsection{Solution insights}
\begin{figure}[t]
    \centering
    \includegraphics[width=1.0\linewidth]{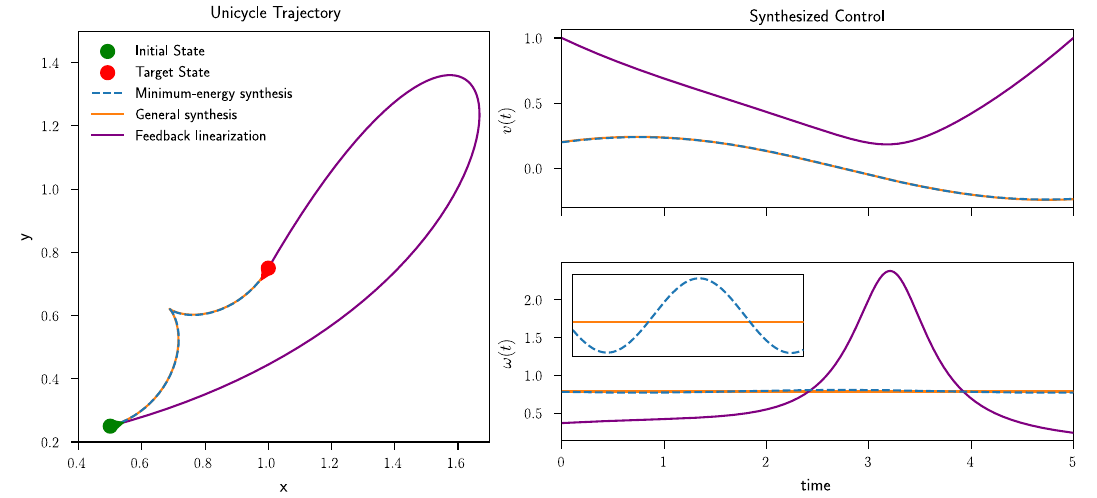}
    \caption{Unicycle Example}
    \label{fig:unicycle}
\end{figure}
\begin{figure}
    \centering
    \includegraphics[width=1.0\linewidth]{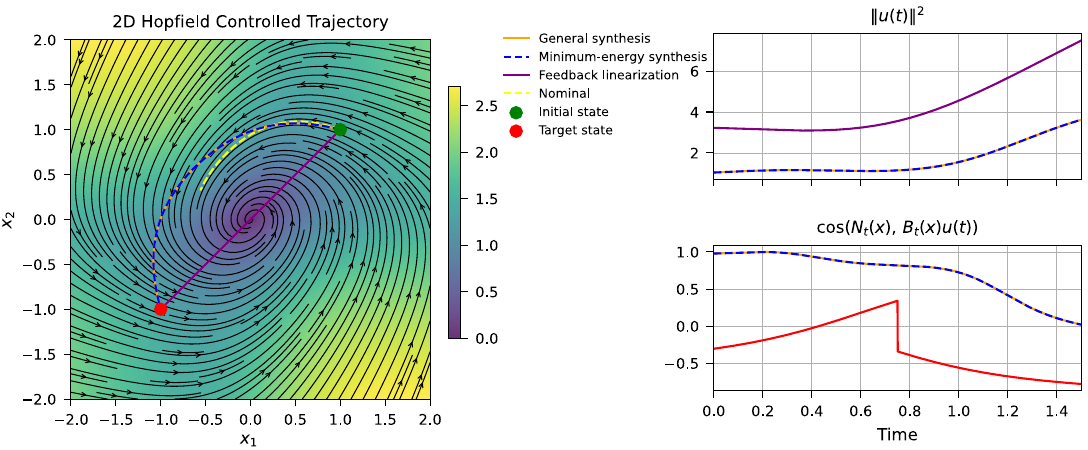}
    \caption{Fully-actuated 2D Hopfield Example}
    \label{fig:hopfield}
\end{figure}
Fig.~\ref{fig:unicycle} compares the synthesized trajectories and control inputs for a unicycle system steered from $(x,y,\theta) = (0.5, 0.25, 
\pi/12)$ to $(1.0, 0.75, 4\pi/3)$ over $[0, 5]$, a non-trivial task requiring simultaneous planar relocation and reorientation. Both general and minimum-energy (ME) syntheses achieve exact 
steering with endpoint errors of $4.17 \times 10^{-12}$ and $9.21 \times 10^{-10}$, respectively. 
Both methods produce visually indistinguishable trajectories (Fig.~\ref{fig:unicycle} left panel). Despite this, the ME synthesis achieves a slightly lower control energy 
($\|u\|^2_2 = 1.80392$) compared to the general synthesis ($\|u\|^2_2 = 1.80416$), as expected. For reference, feedback linearization---based on differential flatness---also achieves exact steering with an endpoint error of $9.14^{-10}$, but consumes 40\% more energy  ($\|u\|^2_2 = 2.555$).
Further examination of the synthesized control inputs (Fig.~\ref{fig:unicycle}, right panel) reveals structurally distinct strategies.
The general synthesis produces a constant $\omega(t)$, while the ME synthesizes a time-varying oscillatory control. 

To further illustrate the energy efficiency of both Gramian syntheses, Fig.~\ref{fig:hopfield} presents results steering a fully-actuated 2D Hopfield system from $(1.0, 1.0)$ to $(-1.0, -1.0)$ over $[0, 1.5]$. For comparison, we use feedback linearization---based on differential flatness---as a naive baseline. All three methods achieve exact steering (error $<10^{-10}$).
However, the control energies differ 
substantially: both Gramian syntheses achieve $\|u\|^2_2 = 1.566$, compared to $\|u\|^2_2 = 2.527$ for feedback linearization, a $38\%$ reduction. The lower energy of both syntheses also aligns with the theory from~\cite{tamekue2025controlanalysis,tamekue2026minimum}---the control from the ME synthesis can coincide with that from the general synthesis. The left panel of Fig.~\ref{fig:hopfield} shows that both Gramian syntheses closely follow the natural flow of the vector field, while feedback linearization takes a more direct but energetically costly path. The upper right panel confirms this, as both Gramian methods maintain significantly lower $\|u(t)\|^2$ throughout. The lower right panel provides an intuitive illustration of the underlying geometric mechanism: both Gramian syntheses maintain high cosine similarity between $B_t(x)u(t)$ and $N_t(x)$, indicating that the synthesized control aligns better with the natural drift. 
\subsection{Controlling Underactuated 100D MINDy via Minimunm Energy Synthesis}  \label{sec:underactuated}
\begin{figure}[t]
    \centering
    \includegraphics[width=1.0\linewidth]{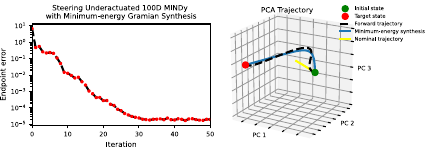}
    \caption{Left: Convergence of 100D mindy with 50D input Right: PCA plot of synthesized trajectory}
    \label{fig:50d}
\end{figure}
We additionally evaluate the synthesis on an underactuated 100D MINDy model with $k = 50$ input channels, a challenging example. 
The initial state was randomly sampled from $\mathcal{N}(\mathbf{0}_{100}, \mathbf{I}_{100})$.
To ensure reachability, the target state was generated by forward simulation with a reference control $u_{\text{sim}}$, constructed as a fifth-order Chebyshev polynomial with coefficients drawn independently from $\mathcal{N}(0, 0.2^2)$. 
The simulation is conducted in double precision with $rtol = 10^{-7}, atol = 10^{-7}$ over 50 Picard iterations.

Fig. \ref{fig:50d} shows the convergence and synthesized control for this underactuated system. The endpoint error decreases steadily, reaching approximately $10^{-5}$ by iteration 50, as expected, slower than the fully actuated case. The synthesized control achieves a lower energy ($\|u(t)\|_2 = 0.727$) compared to the reference control ($\|u_{\text{sim}}(t)\|_2 = 1.103$), a 34\% reduction. 
The PCA projection (Fig. \ref{fig:50d} right) of the trajectory confirms that the synthesized path departs significantly from both the nominal and forward trajectories before converging to the target, reflecting a distinct steering strategy under actuation constraints.

\subsection{ Runtime Scalability of General Synthesis}
We further evaluate the computational scalability of the general synthesis up to 128 dimensions. The scalability here refers to computational 
feasibility in high dimensions rather than characterizing the asymptotic growth rate. 

We construct submodels of dimension $d \in \{2, 8, 32, 64, 128\}$ by randomly selecting $d$ nodes from a 200-dimensional MINDy model estimated from resting-state fMRI data~\cite{chenDynamicalModelsReveal2025a}, preserving the corresponding rows and columns of the synaptic weight matrix $W$, decay matrix $D$, and activation parameters $\alpha$ and $\beta$. The system here is fully controllable with $B_t(x_t)=\idty$. All experiments are conducted on an NVIDIA V100 GPU under double precision with $\text{rtol}=10^{-8}$, $\text{atol}=10^{-10}$, 5000 Simpson quadrature points, and Cholesky-based Gramian inversion. For numerical stability, we regularize the Gramian before inversion by adding $\epsilon\idty$ with $\epsilon={10}^{-6}$. For each model size, we run 5 independent trials. All the trials use the same solver configuration, time horizon $[0, 10]$ and initial state $x^0 = \mathbf{0}_m$. We justify this initial state choice analytically as an equilibrium of the MINDy drift dynamics ($N_t(0) = 0$). The target state $x^1$ is sampled uniformly at random using a model-size- and trial-specific random seed \texttt{jax.random.key} to maximize the result's generalizability.


Fig.~\ref{fig:scale} reports the amortized per-iteration wall-clock time as a function of $d$ for $N_{\max} \in \{1, 2, 5, 10\}$. 
The general synthesis remains computationally feasible at $d=128$ in double precision, with an average amortized per-iteration cost close to 25 seconds, establishing practical tractability at dimensionality relevant for the types of neuroscience models discussed above.
The vertical separation between curves reflects the one-time JIT compilation overhead amortized across iterations, which becomes negligible at large $d$ as per-iteration execution cost dominates.
It is worth noting that, in addition to runtime reduction, JIT also effectively reduced peak memory usage, which is expected behavior for operation optimization~\cite{jax2018github}.

We note that the reported times should be interpreted as evidence of computational feasibility specific to the chosen setup rather than as a general runtime guarantee across all systems and boundary conditions. 
The models here don't have super strong nonlinearity. Systems with strong nonlinear dynamics require significantly more adaptive solver steps, increasing per-iteration cost independently of dimension. Additionally, GPU memory presents a significant yet latent bottleneck not captured by wall-clock time alone.
Thus, we focus on feasibility for large-scale systems rather than on observing a scaling law.
\begin{figure}[t]
    \centering
    \includegraphics[width=1.0\linewidth]{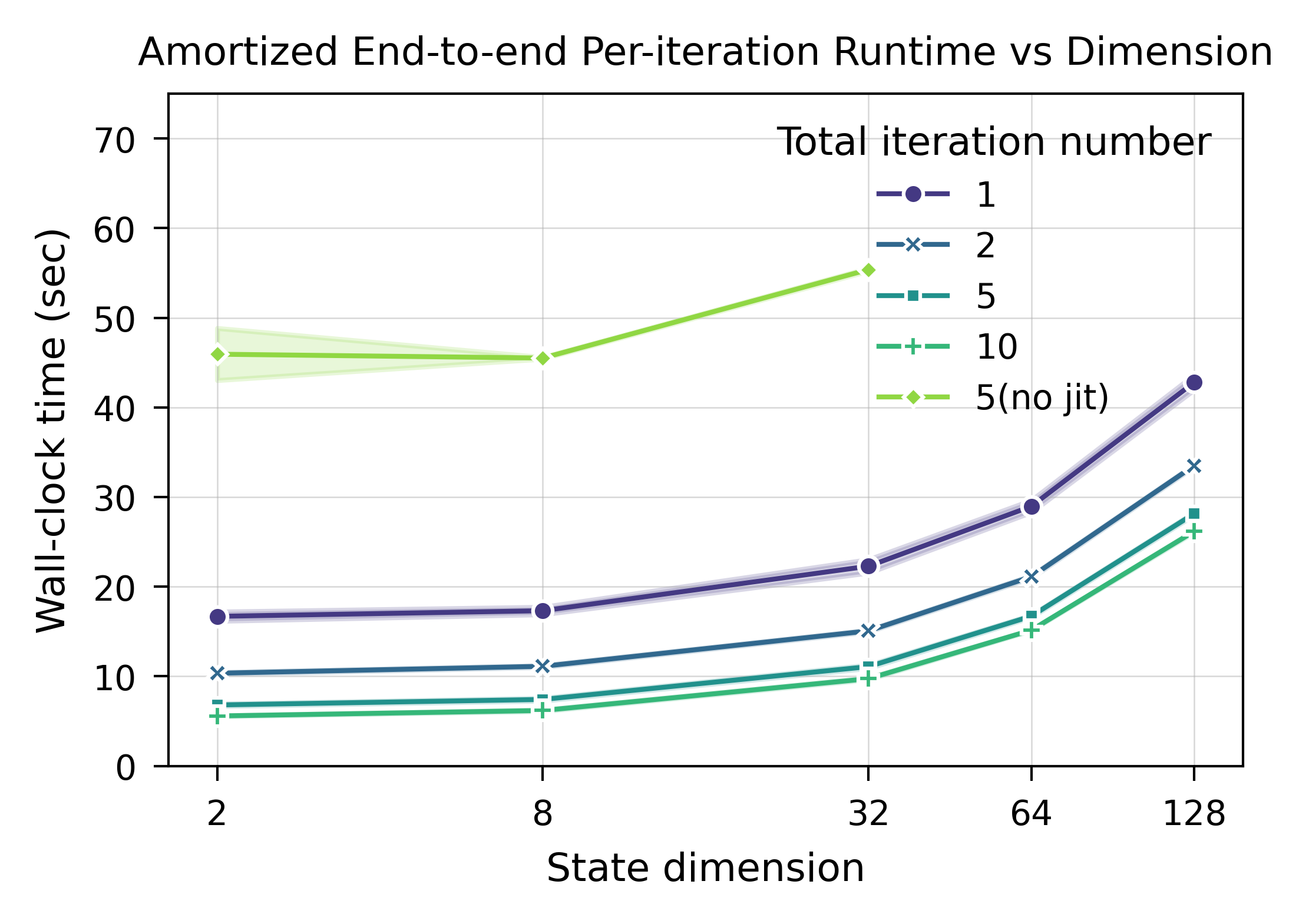}
    \caption{Runtime benchmark of General synthesis. The x-axis is shown on a $\log_2$ scale. 
    Shaded bands indicate standard deviation over 5 runs.}
    \label{fig:scale}
\end{figure}

\section{Conclusions}
This work presents a scalable numerical implementation of the nonlinear Gramian synthesis framework~\cite{tamekue2025controlanalysis,tamekue2026minimum}, addressing the computational bottlenecks inherent to nonlinear synthesis and scaling the method to fully- and under-actuated 100-dimensional brain models estimated from human fMRI data. We provide the first numerical implementation of the minimum-energy Gramian synthesis, documenting the computational foundation underlying~\cite{tamekue2026minimum}. 
Our results confirm that convergence is governed by intrinsic controllability rather than state-space dimension.



\bibliographystyle{ieeetr}
\bibliography{picard,brain_modeling}

\end{document}